\documentclass{gtart_h}

\def\ifplaintex{\expandafter\ifx\csname documentclass\endcsname\relax}

\def\gtp{{\mathsurround=0pt\it $\cal G\mskip-2mu$eometry \&\ 
$\cal T\!\!$opology $\cal P\!$ublications}}  

\def\Addressesr{\bigskip
{\small \parskip 0pt \leftskip 0pt \rightskip 0pt plus 1fil \def\\{\par}
\sl\theaddress\par
\medskip
\rm Email:\stdspace\tt\theemail\hfill\rm Received:\qua\receiveddate \par}}

\def\recd{{\small Received:\qua\receiveddate\ifx\reviseddate\relax
\else\qquad Revised:\qua\reviseddate\fi\par}} 

\def\lognumber#1{\def\thelognumber{#1}}
\def\volumenumber#1{\def\thevolumenumber{#1}}
\def\volumeyear#1{\def\thevolumeyear{#1}}
\def\papernumber#1{\def\thepapernumber{#1}}
\def\pagenumbers#1#2{\def\startpage{#1}\def\finishpage{#2}}
\def\published#1{\def\publishdate{#1}}

\def\received#1{\def\receiveddate{#1}}

\def\accepted#1{\def\accepteddate{#1}}

\long\def\asciiabstract#1{\long\def\theasciiabstract{#1}}

\let\\\par\let\thelognumber\relax\let\thevolumenumber\relax
\let\thepapernumber\relax\let\thevolumeyear\relax\let\startpage\relax
\let\finishpage\relax\let\publishdate\relax\let\receiveddate\relax
\let\reviseddate\relax\let\accepteddate\relax\let\theasciititle\relax
\let\theasciiauthors\relax
\let\theasciiabstract\relax

\let\theasciiemail\relax

\ifplaintex
\font\logobig=cmssbx10 scaled 3836
\font\logomed=cmssbx10 scaled 2557
\else
\font\logobig=cmssbx10 scaled 4200
\font\logomed=cmssbx10 scaled 2800
\fi

\long\def\makeagttitle{   
\count0=\startpage
\agt\hfill      
\hbox to 45truept{\vbox to 0pt{\vglue -13truept{\logomed A\kern -.37em{\logobig 
T}\kern -.38em G}\vss}\hss}
\break
{\small Volume \thevolumenumber\ (\thevolumeyear)
\startpage--\finishpage\nl
Published: \publishdate}

\vglue .25truein

{\parskip=0pt\leftskip 0pt plus
1fil\def\\{\par\smallskip}{\Large\bf\thetitle}\par\medskip} \vglue
0.05truein

{\parskip=0pt\leftskip 0pt plus 1fil\def\\{\par}{\sc\theauthors}
\par\medskip}%
 
\vglue 0.03truein 

{\small\leftskip 25truept\rightskip 25truept{\bf Abstract}\stdspace\theabstract

{\bf AMS Classification}\stdspace\theprimaryclass
\ifx\thesecondaryclass\relax\else; \thesecondaryclass\fi\par
{\bf Keywords}\stdspace \thekeywords\par}\vglue 7truept

}   

\ifplaintex
\font\phead=cmsl9 scaled 950
\font\pnum=cmbx10 scaled 913
\font\pfoot=cmsl9 scaled 950
\headline{\vbox to 0pt{\vskip -4.5mm\line{\small\phead\ifnum
\count0=\startpage ISSN 1472-2739 (on-line) 1472-2747 (printed)
\hfill {\pnum\folio}\else\ifodd\count0\def\\{ }%
\ifx\theshorttitle\relax\thetitle\else\theshorttitle\fi\hfill{\pnum\folio}
\else\def\\{ and }{\pnum\folio}\hfill\ifx\theshortauthors\relax\theauthors
\else\theshortauthors\fi\fi\fi}\vss}}
\footline{\vbox to 0pt{\vglue 0mm\line{\small\pfoot\ifnum\count0=\startpage
\copyright\ \gtp\hfill\else
\agt, Volume \thevolumenumber\ (\thevolumeyear)\hfill\fi}\vss}}
\else
\font=cmsl9 scaled 1050
\font\lnum=cmbx10 
\font=cmsl9 scaled 1050
\makeatletter
\def\@oddhead{{\small\lhead\ifnum\count0=\startpage ISSN 1472-2739 
(on-line) 1472-2747 (printed)\hfill {\lnum\number\count0}\else\ifodd\count0
\def\\{ }\ifx\theshorttitle\relax \thetitle \else\theshorttitle\fi\hfill
{\lnum\number\count0}\else\def\\{ and }{\lnum\number\count0}
\hfill\ifx\theshortauthors\relax 
\theauthors\else\theshortauthors\fi\fi\fi}}\def\@evenhead{\@oddhead}
\def\@oddfoot{\small\lfoot\ifnum\count0=\startpage\copyright\ \gtp\hfill\else
\agt, Volume \thevolumenumber\ (\thevolumeyear)\hfill\fi}
\def\@evenfoot{\@oddfoot}
\makeatother
\fi
\let\maketitlepage\makeagttitle
\let\makeshorttitle\maketitlepage
\let\maketitle\maketitlepage

\newwrite\gtoutfile
\long\gdef\makeheadfile{  
{\def\\{, }\def\s{ }
\immediate\openout\gtoutfile head.xxx
\immediate\write\gtoutfile{Proxy-for: \ifx\theasciiauthors\relax
\theauthors\else\theasciiauthors\fi\s<\ifx\theasciiemail\relax\theemail\else\theasciiemail\fi>}
\immediate\write\gtoutfile{\noexpand\\}
\immediate\write\gtoutfile{Authors: \ifx\theasciiauthors\relax
\theauthors\else\theasciiauthors\fi}
{\def\\{ }\immediate\write\gtoutfile{Title: \ifx\theasciititle\relax
\thetitle\else\theasciititle\fi}}
\immediate\write\gtoutfile{Subj-class: GT or SG, GR etc}
\immediate\write\gtoutfile{MSC-class: \theprimaryclass\ifx\thesecondaryclass\relax\else, \thesecondaryclass\fi}
\immediate\write\gtoutfile{Journal-ref: Algebr. Geom. Topol. \thevolumenumber\s
(\thevolumeyear) \startpage-\finishpage}
\immediate\write\gtoutfile{Comments: Published by Algebraic and
Geometric Topology at}
\immediate\write\gtoutfile{\s\s\s  http://www.maths.warwick.ac.uk/agt/AGTVol\thevolumenumber/agt-\thevolumenumber-\thepapernumber.abs.html}
\immediate\write\gtoutfile{\noexpand\\}
\immediate\write\gtoutfile{}
\ifx\theasciiabstract\relax
\immediate\write\gtoutfile{\theabstract}\else
\immediate\write\gtoutfile{\theasciiabstract}\fi
\immediate\write\gtoutfile{}
\immediate\write\gtoutfile{\noexpand\\}
\immediate\write\gtoutfile{}
\immediate\closeout\gtoutfile}}  

\def\maketitlepage{\makeagttitle\makeheadfile}
\let\makeshorttitle\maketitlepage
\let\maketitle\maketitlepage

\lognumber{31}
\volumenumber{5}
\volumeyear{2005}
\papernumber{31}
\pagenumbers{741}{750}
\received{16 November 2004} 
\accepted{6 March 2005}
\published{23 July 2005}

\usepackage{amssymb} 
\usepackage{graphicx}

\newtheorem{theorem}{Theorem}

\newtheorem{prop}{Proposition}

\newtheorem{ques}{Question}

\theoremstyle{definition}
\newtheorem*{rem}{Remark}

\newcommand{\B}{\mathcal B}

\newcommand{\C}{\mathbb C}
\newcommand{\CP}{\mathbb{CP}^1}

\newcommand{\p}{\partial}

\newcommand{\pM}{\partial M}
\newcommand{\pn}{\pi}

\newcommand{\Q}{\mathbb Q}
\newcommand{\R}{\mathbb R}

\newcommand{\SLC}{\mbox{SL}_2(\C)}

\newcommand{\tXo}{\widetilde{X}_0}

\newcommand{\Z}{\mathbb Z}

\begin{document}

\title{Boundary slopes (nearly) bound cyclic slopes}

\author{Thomas W. Mattman}
\address{Department of Mathematics and Statistics, California State
University, Chico\\Chico, CA95929-0525, USA}
\email{TMattman@CSUChico.edu}

\begin{abstract}  
Let $r_{m}$ and $r_{M}$ be the least and greatest finite boundary
slopes of a hyperbolic knot $K$ in $S^3$. We show that any cyclic
surgery slopes of $K$ must lie in the interval $( r_{m} - 1/2, r_{M} +
1/2 )$.
\end{abstract}

\asciiabstract{%
Let r_m and r_M be the least and greatest finite boundary slopes of a
hyperbolic knot K in S^3. We show that any cyclic surgery slopes of K
must lie in the interval (r_m - 1/2, r_M + 1/2).}


\maketitle

\section{Introduction}
In \cite{M1,M2} we observed that the Seifert surgeries of 
$(-3,3,n)$ pretzel knots follow an interesting pattern, summarised in the
table below. For each positive integer $n$, the Seifert surgeries of the $(-3,3,n)$ pretzel
lie between the boundary slopes $0$ and $8/(n+1)$. Indeed, all integral slopes 
in the interval $(0, 8/(n+1))$ are Seifert.
\def\Strut{\vrule width 0pt height 12pt depth 7pt}
\begin{center}
\begin{tabular}{|c|c|c|c|c|c|c|c|} \hline
\Strut$n$ & $1$ & $2$ & $3$ & $4$ & $5$ & $6$ & $\geq 7$\\ \hline
\Strut$8/(n+1)$ & $4$ & $8/3$ & $2$ & $8/5$ & $4/3$ & $8/7$ & $\leq 1$ \\ \hline
\Strut Seifert Surgeries & $1,2,3$ & $1,2$ & $1$ & $1$ & $1$ & $1$ & none \\ \hline
\end{tabular}
\end{center}
Some other famous knots also share this pattern; all integral slopes 
between two boundary slopes (shown in bold typeface) are Seifert:
\begin{center}
\begin{tabular}{|l|c|} \hline
\Strut Knot & Non-trivial exceptional surgeries \\
\hline
\Strut Figure 8 & $\mathbf{-4}$, $-3$, $-2$, $-1$, $\mathbf{0}$, $1$,
$2$, $3$, $\mathbf{4}$ \\
\Strut(-2,3,7) & $\mathbf{16}$, $17$, $18$, $\mathbf{37/2}$, $19$,
$\mathbf{20}$ \\
\Strut Twist Knots & $\mathbf{0}$, $1$, $2$, $3$, $\mathbf{4}$ \\
 \hline
\end{tabular}
\end{center}
Based on this and other evidence, Kimihiko Motegi posed the question:
``Are Seifert surgeries bounded by boundary slopes?'' In other words, if
slope $t$ is a Seifert surgery slope, are there necessarily boundary slopes 
$r_m$ and $r_M$ with $r_m \leq t \leq r_M$?  
In \cite{IMS}, we construct a parameter $c$, such that 
$r_m - c \leq t \leq r_M + c$. More
precisely,  we can reformulate \cite[Corollary 3]{IMS} as follows. (In part 3 of
Theorem~\ref{thC3IMS}, the constant
$s$ is the minimal Culler-Shalen norm defined in Section 2 below while $A$ counts the 
characters of non-abelian representations of the knot exterior that factor through the
surgery. We refer the reader to \cite{IMS} for details.)

\begin{theorem}[Corollary 3 of~\cite{IMS}]\label{thC3IMS}
Let $r_m$ and $r_M$ be the least and greatest finite boundary slopes of a hyperbolic knot
$K$ and
$t$ a non-trivial exceptional surgery slope. Then $r_m - c \leq t \leq r_M + c$ where $c$
depends on the type of slope $t$.  
\begin{enumerate}
\item If $t$ is cyclic, $c = 1$.
\item If $t = a/b$ is finite, $c = 3/b$.
\item If $t = a/b$ is a Seifert fibred slope, $c = (1+{2A}/s)/b$.
\end{enumerate}
\end{theorem}

In this formulation, Motegi's conjecture corresponds to showing $c = 0$ for a Seifert fibred 
surgery. That $c = 1$ for a cyclic surgery was first shown by Dunfield~\cite{Du}.
In the current article, we show that for a cyclic surgery, we can take $c = 1/2$.

\begin{theorem} \label{thmain}
If $t$ is a non-trivial cyclic surgery on a hyperbolic knot $K$ in $S^3$ and $r_m$ and $r_M$
are the least and greatest finite boundary slopes of $K$, then
$r_m - \frac12 < t < r_M + \frac12$.
\end{theorem}

Moreover, the theorem applies more generally to a hyperbolic knot in a manifold with
cyclic fundamental group whose exterior satisfies $H^1(M, \Z_2) = \Z_2$. 

In light of Theorem~\ref{thC3IMS}, it is natural to extend Motegi's question about Seifert
surgeries to exceptional surgeries in general:

\begin{ques} Do boundary slopes bound exceptional slopes? \end{ques}

In other words, for a hyperbolic knot $K$ in $S^3$ do all non-trivial exceptional surgery
slopes lie in the interval ${[}r_m, r_M{]}$ between the least and greatest finite boundary
slopes?  Note that, by the Cabling Conjecture, $K$ should have no reducible surgeries and
that toroidal surgeries are themselves boundary slopes and 
will, therefore, necessarily lie in ${[} r_m, r_M {]}$. The real question is whether other
types of exceptional  surgeries (i.e., cyclic, finite, Seifert fibred) must also lie in this
interval. It is the Seifert case that gives the largest values for $c$ in
Theorem~\ref{thC3IMS} and, in any case, cyclic and finite surgeries are
thought to be examples of Seifert surgeries. Thus, an affirmative answer
to Motegi's question  is likely to imply the same for all types of
exceptional surgeries. 


In Section 2 we provide definitions and discuss the geometry of the Culler-Shalen norm.
In Section 3, we prove Theorem~\ref{thmain}.

\section{Definition, geometry of the Culler-Shalen norm}

Let $K$ be a hyperbolic knot in $S^3$ and let $M = S^3 \setminus N(K)$ denote the knot
exterior. Fixing the usual meridian, longitude basis $\{ \mu, \lambda \}$, the element
$\gamma = a
\mu + b
\lambda$ of $H_1(\pM ; \Z)$ will be represented as $(a,b)$. This class can be identified with
the {\em ``slope"} $r_{\gamma} = a/b$ in 
$\Q \cup \{ \frac10 \}$. We will occasionally wish to change our framing which amounts
to replacing $\lambda$ by $k \mu  + \lambda$ and to changing coordinates
by
$(a,b) \mapsto (a-bk,b)$.

Let $M(r)$ denote the manifold obtained by Dehn surgery along
slope $r$ (i.e., $M(r)$ is constructed by attaching a solid torus to $M$ such that the
boundaries of meridional disks are curves of slope $r$ in $\pM$). We will call
$r$ a {\em cyclic} (respectively {\em finite}) {\em slope} if
$\pi_1(M(r))$ is cyclic (resp. finite). If $M(r)$ admits the structure of a Seifert fibred
space, we call $r$ a  {\em Seifert fibred slope}. Since $M(\frac10) = S^3$, we refer to
meridional surgery along slope $r_\mu = \frac10$ as {\em trivial surgery}.

If there is an essential surface $\Sigma$ in $M$ that meets $\partial M$ in a
non-empty set of parallel curves of slope $r$, we call $r$ a {\em boundary slope}. If there
is such a $\Sigma$ that is not a fibre in a fibration of $M$ over $S^1$, $r$ is a {\em strict
boundary slope}. For example, by applying the loop theorem to a Seifert surface of $K$, we
observe that $0$ is a boundary slope. We will say $r$ is a {\em finite boundary slope} if it
is a boundary slope and $r \neq \frac10$.

The proof of Theorem~\ref{thmain} depends on the geometry of the Culler-Shalen norm of $K$.
We introduce some of the main properties of this norm and refer the reader to
\cite[Chapter 1]{CGLS} for a more complete account.

Let $R = \mbox{Hom}(\pn, \SLC)$ denote the set of
$\SLC$-representations of the fundamental group $\pn$ of $M$.
Then $R$ is an affine algebraic set, as is $X$, the set of characters
of representations in $R$.

For $\gamma \in \pn$, define the regular function
$I_{\gamma}:X \to \C$ by $I_{\gamma}(\chi_{\rho}) = \chi_{\rho}(\gamma) =
\mbox{trace}(\rho(\gamma))$.
By the Hurewicz isomorphism,
a class $\gamma \in L =  H_1(\pM, {\Z})$ determines an element
of $\pi_1(\pM)$, and therefore an element of $\pn$ well-defined
up to conjugacy. A {\em norm curve} is a one-dimensional irreducible component of $X$
on which no $I_\gamma$ ($\gamma \in L \setminus \{ 0 \}$) is constant. For example, the
irreducible component, $X_0$, that contains the character of the holonomy representation is
a norm curve.

The terminology reflects the fact that we may associate to $X_0$ a norm $\| \cdot \|$ on
$H_1( \pM, \R)$ called a {\em Culler-Shalen norm} in the following manner. Let $\tXo$ be the
smooth projective model of $X_0$ which is birationally equivalent to $X_0$. The 
birational map is regular at all but a finite number of points of $\tXo$ which are called
{\em ideal points} of $\tXo$.  The function
$f_{\gamma} = I_{\gamma}^2 -4$ is again regular and so can be pulled back to
$\tXo$. For $\gamma \in L$, $\| \gamma \|$ is the degree of
$f_{\gamma} : \tXo \to \CP$.  The norm is extended to $H_1(\pM, \R)$ 
by linearity. 

Let $s = \min_{0 \neq \gamma \in H_1 (\pM, \Z)}  \| \gamma \| $ denote the 
{\em minimal norm}.
The norm disc of radius $s$ is a convex, finite-sided polygon $P$ that is symmetric about
the origin. We will call $P$ the {\em fundamental polygon}. The ideal points of $\tXo$ can
be associated with a set $\B$ of strict boundary slopes of the knot and the vertices of $P$
occur at rational multiples of the classes of slopes in $\B$. It follows that $\B$ must
contain at least two slopes.  One of the main results of \cite{CGLS} is that if $r_{\gamma}$
is a cyclic slope that is not a strict boundary slope then $\| \gamma \| = s$. Moreover,
$r_\gamma$ is either integral or trivial (i.e., $r_\gamma = \frac10$).

\section{Proof of Theorem~\ref{thmain}}

We will prove two propositions before coming to the proof of the theorem. 
The main external inputs for our argument are the Cyclic Surgery Theorem~\cite{CGLS} and
Theorems 4.1 and 4.2 of \cite{Du}. Thus, although we formulate our results in terms of a
knot in $S^3$, they carry over to the case of a hyperbolic knot in a manifold with
cyclic fundamental group whose exterior satisfies $H^1(M, \Z_2) = \Z_2$. 

\begin{prop} 
Let $K$ be a hyperbolic knot in $S^3$. Let $r_M$ be the greatest finite boundary slope of
$K$. Suppose $\gamma$ is a non-trivial cyclic class with $r_\gamma = n > 0$.
Then $n \leq r_M + \frac12$.
\end{prop} 

\begin{rem}
In fact, we will show that $n \leq r + \frac12$ where $r$ is the
greatest finite  boundary slope associated to the norm curve $X_0$. In
particular,
$r$ is a strict boundary slope.
\end{rem}

\proof Let $X_0$ be the norm curve that contains the character of the holonomy
representation and
let $\B$ be the associated set of  boundary slopes. If $r_\gamma \in \B$, then
$n = r_\gamma \leq r_M$. So we may assume $r_\gamma \not\in \B$.

\begin{figure}[ht!]
\begin{center}
\includegraphics[width=0.7\hsize]{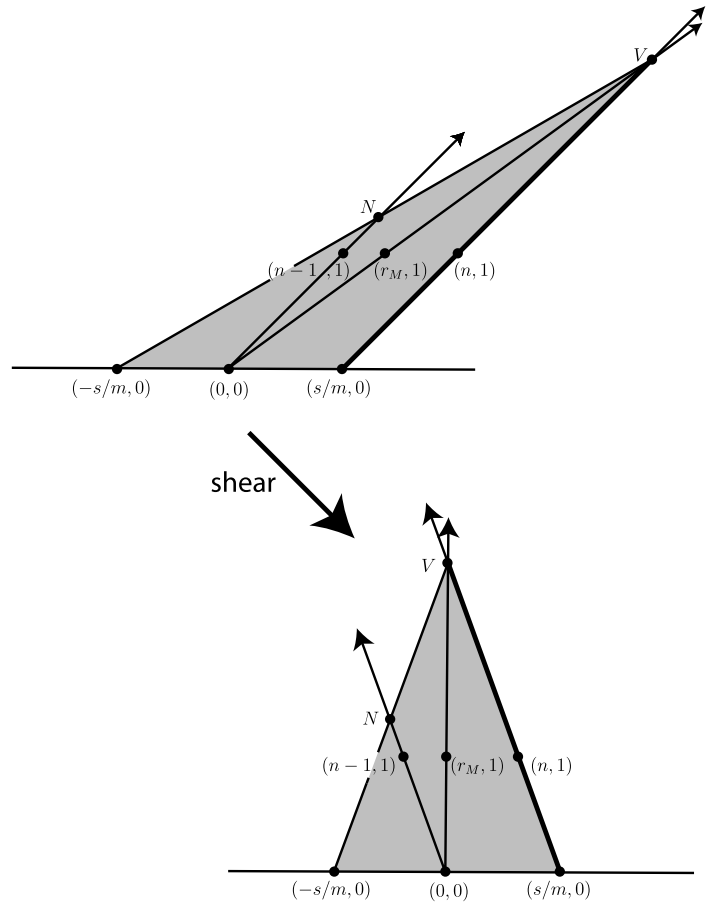}
\caption{\label{figprop1}%
The geometry of $P$ assuming $n > r_M + \frac12$}
\end{center}
\end{figure}

Suppose, for a contradiction, that $n > r_M + \frac12$.
Without loss of generality, we may assume $r_M \in \B$ (otherwise replace
$r_M$ by the greatest finite boundary slope in $\B$).
Our goal is to argue that $(n-1 ,1)$ is in the interior of $P$ (see Figure~\ref{figprop1})  
where $P$ is the fundamental polygon of the Culler-Shalen norm $\| \cdot \|$ associated
to $X_0$.
Let us outline the argument.
First, it has already been shown that $r_M$ must be non-integral~\cite{Du}.
Construct the line from $(s/m,0)$ through $(n,1)$, where
$m = \| \mu \|$ and $s$ is the minimal norm.
This line will form part of the boundary of $P$ and continues to the vertex
$V$ determined by the boundary slope
$r_M$. By convexity, the line joining $(-s/m,0)$ to $V$ is also in $P$. It follows that
$(n-1,1)$ is in the interior of $P$ 
so that $\|(n-1,1) \| < s$. (This is ``obvious'' if one applies a shear as in 
Figure~\ref{figprop1}.)
This is a contradiction as $s$ is defined to
be the minimal norm. Thus, we conclude $n \leq r_M + \frac12$.

Now let's fill in the details.
Since $n > r_M + \frac12$, if $r_M$ is an integer, then $r_M \leq n-1$. However, 
by \cite[Theorem 4.2]{Du} there is a strict boundary slope $r_\delta \in \B$
with $|n - r_{\delta}|<1$. Then $r_M < r_{\delta}$ in
contradiction to the choice of $r_M$. Therefore,
$r_M$ is not an integer.  Moreover, since there is a strict boundary
slope $\delta$ with $|n-r_{\delta}| < 1$ we
have $n-1 < r_M < n$. 

We next construct the vertex $V$ of $P$ corresponding to the boundary slope 
$r_M$. 
Let $\| \mu \| = m$.
(If $r_\mu = \frac10 \not\in \B$, then $m = s$.) 
The
point $(s/m,0)$ is then in $\p P$.
Since $r_M$ is maximal among
slopes in $\B$, and $r_M < n < \infty$, the segment joining
$(s/m,0)$ and $(n,1)$ is part of the boundary
of $P$. This segment has equation $y = (x-s/m)/(n-s/m)$. It continues to
the line $y = x/r_M$. (Since there are
no strict boundary slopes between $n$ and $r_M$, the segment has no
vertex before it reaches the line $y =
x/r_M$ corresponding to the boundary slope $r_M$.) These lines
meet at the point $$V = \frac{s/m}{r_M + s/m -
n}(r_M,1)$$ which is therefore a vertex of $P$. 

Since $P$ is
convex, the segment joining $V$ and $(-s/m,0)$
(both in $P$) is contained in $P$. We argue that the
point $N$ where this segment crosses $y = x/(n-1)$ is above the
line $y=1$.
Indeed, the segment has the equation $$y =
\frac{x + s/m}{2 r_M + s/m - n}.$$ It
meets the line $y = x/(n-1)$ at the point $$N = \frac{s/m}{2(r_M -
n) + 1 + s/m} (n-1,1)$$ which is therefore in
$P$. Let $y_N$ denote the $y$ coordinate of $N$. 
\begin{eqnarray*} n > r_M + \frac12 &
\Rightarrow & 0 > 2(r_M - n) +1 \\ 
& \Rightarrow & s/m + 2(r_M - n) + 1 < s/m \\
& \Rightarrow & y_N = \frac{s/m}{2(r_M -n) + 1 + s/m} > 1
\end{eqnarray*} 
Since $y_N > 1$, 
the point $(n-1,1)$ is in the interior of $P$
and, therefore, 
\mbox{$\| (n-1,1) \|$} $< s$. This is a contradiction as $s$ is defined to be the
minimal norm. We conclude that $n \leq r_M +
\frac12$. \endproof

We will now show that 
Proposition~1 can be strengthened to a strict inequality
if the meridian is not a strict boundary class for the norm curve $X_0$. The
argument makes use of the idea of the diameter
$D$ of the set of boundary slopes. Culler and Shalen~\cite{CS3} showed that $D \geq 2$.
This inequality is sharp by an example of Dunfield of a knot in a manifold
with cyclic fundamental group (see~\cite{CS3}). For hyperbolic knots in
$S^3$, the smallest known diameter is $D = 8$ for
the Figure 8 knot.

\begin{prop} 
Let $X_0$ be the norm curve containing the character of the holonomy representation
for the hyperbolic knot $K \in S^3$. Let $\B$ be the associated boundary
slopes and suppose that $r_\mu = \frac10 \not\in
\B$. Let $\gamma$ be a non-trivial cyclic class with $r_\gamma = n > 0$.
Let $r_m$ and $r_M$ denote the least and greatest boundary slopes in 
$\B$. Let $D = \mbox{Diam}(\B) = r_M - r_m$. Then
$n \leq r_M + 1 - \frac12(D - \sqrt{D(D-2)}) < r_M +
\frac12$.
\end{prop}

\begin{rem}
The difference between $n$ and $r_M$ goes to zero as $D$ approaches
$2$. For
$D = 8$ (Figure 8 knot) we have $1-\frac12(D - \sqrt{D(D-2)}) = -3 + 2\sqrt{3} \approx
0.46$.
\end{rem}

\begin{figure}[ht!]
\begin{center}
\includegraphics[width= 0.7\hsize]{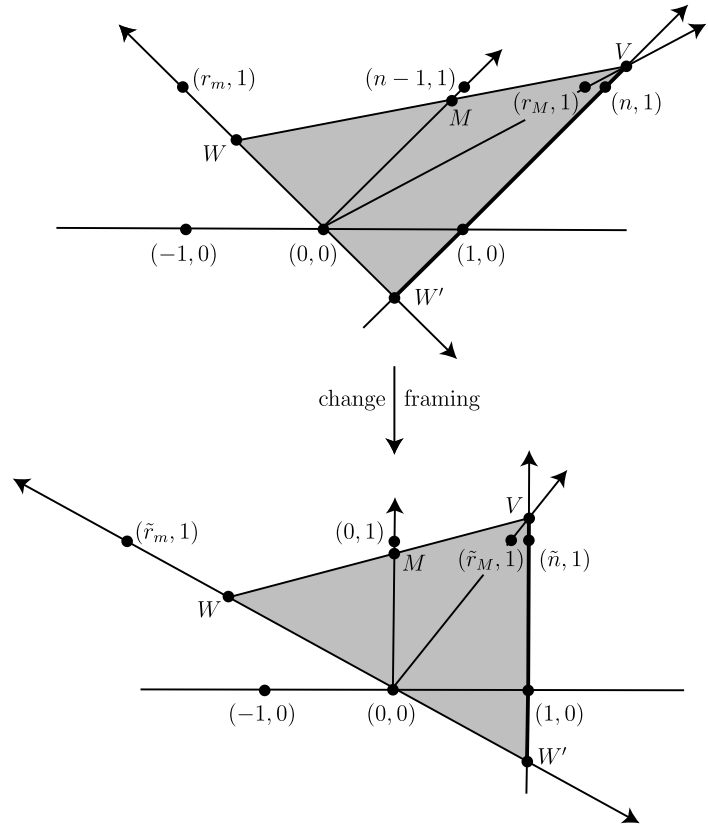}
\caption{\label{figprop2}%
The geometry of $P$ assuming $\| \mu \| = s$}
\end{center}
\end{figure}

\proof Now $\| \mu \| = s$ and $\pm (1,0) \in \p P$. 
If $n \leq r_M$, the proposition holds. So, using Proposition~1, we'll
assume $r_M < n \leq r_M + \frac12$.
Let's change the framing so that $n$ becomes $1$. We will use $\sim$ to
refer to measurements in the new framing. Thus, $\tilde{n} = 1$,
$\tilde{r}_M = r_M - n + 1$ and $\tilde{r}_m = r_m - n + 1$.
The line
through $(\tilde{n},1)$ and $(1,0)$ is then vertical and 
and meets $y = x/{\tilde{r}_M}$ at $V = (1, 1/{\tilde{r}_M})$ (see Figure~\ref{figprop2}).

The segment in $\p P$ which passes through $V$, $(\tilde{n},1)$, and $(1,0)$
continues to  the line $y = x/{\tilde{r}_m}$ as there are no boundary slopes
between $\tilde{r}_m$ and $\tilde{r}_M$ to provide a vertex. The intersection point
$W' = (1,1/ \tilde{r}_m)$
is therefore a vertex of $P$ as is its reflection $W = -W'$.

Since $P$ is convex, the segment joining $V$ and $W$ is contained in $P$.
It meets the line $y = 0$ at the point 
$$M = (0, \frac{\tilde{r}_m - \tilde{r}_M}{2 \tilde{r}_M \tilde{r}_m }).$$
Since $\|(0,1)\| \geq s$, the $y$ coordinate of $M$ cannot exceed $1$ (note
that $D = \tilde{r}_M - \tilde{r}_m$):
$$ 1 \geq \frac{D}{2 \tilde{r}_M (D - \tilde{r}_M )}.$$ 
By the previous proposition, 
$\tilde{r}_M \geq \frac12$. Since $n \geq r_M$, we have 
$\tilde{r}_M \leq 1$. Recall \cite{CS3} that $D \geq 2$. 
\begin{eqnarray*}
1 \geq \frac{D}{2 \tilde{r}_M (D - \tilde{r}_M)} & \Rightarrow &
2 \tilde{r}_M (D - \tilde{r}_M) \geq D \\
& \Rightarrow & 2 \tilde{r}_M^2 -2D \tilde{r}_M + D \leq 0 \\
& \Rightarrow & \tilde{r}_M \geq \frac12 (D - \sqrt{D(D-2)})
\end{eqnarray*}
Since $\tilde{r}_M = r_M - n + 1$, we have, $n \leq r_M + 1 - \frac12(D - \sqrt{D(D-2)})$,
as required.
\endproof

We are now in a position to prove Theorem~\ref{thmain}.

Let $K$ be a hyperbolic knot in $S^3$ and
let $r_\gamma = n$ be a non-trivial cyclic surgery slope. Then, by the Cyclic Surgery
Theorem~\cite{CGLS},
$n$ is an integer and, without loss of generality,  we may assume $n \geq 0$. Since $0$ is a
boundary slope, 
$r_m \leq 0$ so that $r_m - \frac12 < n$. Similarly, if $n = 0$, $r_M + \frac12 > n$.
Thus, in order to prove Theorem~\ref{thmain},
it is enough to show $n < r_M + \frac12$ when $0< n = r_\gamma$ is a cyclic surgery
slope.

By Proposition~1, $n \leq r_M + \frac12$, so it remains only to show that equality is not
possible. Suppose then (for a contradiction) that $n = r_M + \frac12$ and change the framing
so that $n$ goes to $\tilde{n} =1$.
Then $\tilde{r}_M = \frac12$. Following the argument in the proof of 
Proposition~1, the line from the point $E = (s/m,0)$ through $V$ is part of the boundary of
$P$. Since $m = \| \mu \|$, $m$ is at least as big as the minimal norm $s$.  Thus,
$E$ lies on the half open segment $( (0,0), (1,0) {]}$. Similarly, $ \| (1,2) \| \geq s$
and, therefore, $V$, which lies on the boundary of $P$, must be in the half open
segment $((0,0), (1,2) {]}$. However, the line through $E$ and $V$ also passes through  
$(\tilde{n},1) = (1,1)$. The only consistent way to account for all these facts is
to have $V = (1,2)$ and $E = (1,0)$. In other words, $\| \mu \| = s$. Reviewing the
argument of Proposition~2, we see that $\| \mu \| = s$ is exactly the 
extra input needed to deduce that $n < r_M + \frac12$. Thus, we
conclude that $n < r_M +
\frac12$. This is absurd since we began by assuming $n = r_M + \frac12$. The contradiction
shows that, in fact, equality is not possible in Proposition~1. This completes the proof of 
Theorem~\ref{thmain}.

\rk{Acknowledgments}
I am grateful to Kimihiko Motegi for suggesting the question that motivated this study 
and to Masaharu Ishikawa and Koya Shimokawa for helpful conversations. I thank the
referee for many suggestions that significantly improved this paper.

\Addressesr

\end{document}